\documentclass{amsart}
\usepackage[latin1]{inputenc}
\usepackage{latexsym}
\usepackage{amssymb}
\usepackage{amscd}

\def\ent#1{ [ \hspace{-.16 cm}{[} \,#1\, ]\hspace{-.16 cm}{]}}

\def\barr{\begin{array}}
\def\earr{\end{array}}

\def\b{\mathbf{b}}

\newcommand{\CC}{\mathbb{C}}
\newcommand{\RR}{\mathbb{R}}
\newcommand{\NN}{\mathbb{N}}
\newcommand{\MM}{\mathbb{M}}

\def\NN{\mathbb{N}}

\def\bma{\left[\begin{array}}
\def\ema{\end{array}\right]}
\def\ben{\begin{enumerate}}
\def\een{\end{enumerate}}

\newtheorem{fed}{Definition}[section]
\newtheorem{teo}[fed]{Theorem}
\newtheorem*{teo*}{Theorem}

\newtheorem{cor}[fed]{Corollary}
\newtheorem{pro}[fed]{Proposition}
\theoremstyle{definition}
\newtheorem{rem}[fed]{Remark}

\newtheorem{exa}[fed]{Example}

\def\la{\lambda}
\def\al{\alpha}

\def\M{\mathbb{M}}

\def\cF{\mathcal{F}}
\def\cG{\mathcal{G}}

 \DeclareMathOperator{\tr}{tr}

\newcommand{\f}{\cF =\{f_i\}_{i=1}^p}
\newcommand{\g}{\cG =\{g_i\}_{i=1}^r}
\newcommand{\gyf}{\cF \cup \cG }

\newcommand{\anorm}{\mathbf{a}=\{a_i\}_{i\in \NN}}
\def\Aa{\mathbf{a}}

\newcommand{\Bb}{\mathbf{b}}

\newcommand{\suma}[1]{\sum_{i=1}^{#1} a_i}
\newcommand{\fop}[1]{S^{#1}}
\newcommand{\sumla}[1]{\sum_{i=1}^{#1} \la_{n-i+1}}
\newcommand{\comp}[1]{ $\cF$ is  ($\Aa, #1 $)-completable}
\newcommand{\sumal}[1]{\sum_{i=1}^{#1} (a_i +\la_{n-i +1})}

\newcommand{\pint}[1]{\displaystyle \left \langle #1 \right\rangle}

\newcommand{\hil}{\mathcal{H}}

\newcommand{\posop}{L(\mathcal{H})^+}

\begin{document}\baselineskip 5mm

\title[Tight frame completions]{Tight frame completions with prescribed norms.}
\author{P. Massey and  M. Ruiz}
\address{P. Massey and  M. Ruiz, Departamento de Matem\'atica,
Universidad Nacional de La Plata, 50 y 115 (1900), La Plata,
Argentina and IAM-CONICET, Saavedra 15 (1083), Buenos Aires,
Argentina.}

\email{massey@mate.unlp.edu.ar}

\email{mruiz@mate.unlp.edu.ar}

\date{}
\thanks{Partially supported by CONICET (PIP 2083/00), UNLP (11 X350)
and ANPCYT (PICT03-9521)}

\begin{abstract}
Let $\hil$ be a finite dimensional (real or complex) Hilbert space
and let $\{a_i\}_{i=1}^\infty$ be a non-increasing sequence of
positive numbers. Given a finite sequence of vectors $\f$ in
$\hil$ we find necessary and sufficient conditions for the existence of
 $r\in \NN\cup\{\infty\}$ and a Bessel sequence $\g$
in $\hil$ such that $\cF\cup\cG$ is a tight frame for $\hil$ and
$\|g_i\|^2=a_i$ for $1\leq i\leq r$. Moreover, in this case we
compute the minimum $r\in \NN\cup\{\infty\}$ with this
property. Using recent results on the Schur-Horn
theorem, we also obtain a not so optimal but
algorithmic computable (in a finite numbers of steps) tight
completion sequence $\cG$.

\medskip


\medskip

\end{abstract}

 \maketitle

\section{Introduction.}

In recent years, the study of frames in finite dimensional Hilbert spaces has been motivated by a large variety of applications, such as signal processing, multiple antenna
coding, perfect reconstruction filter banks, and Sampling Theory.

Some particular frames, called {\sl tight} frames, are of special
interest since they allow simple reconstruction formulas.
For practical purposes, is often useful to obtain tight
frames with some extra ``structure", for example with the norms
of its elements prescribed (controlled) in advance.

In \cite{FWW} D. Feng, L. Wang and Y. Wang considered the problem
of computing tight completions of a given set of vectors. More
explicitly, given a finite sequence $\f$ of vectors in $\hil$,
 how many vectors we have to add in order to obtain a
tight frame, and how to find those vectors?. Theorem 1.1 in
\cite{FWW} provides a complete answer to this question. But when
the norms of the additional vectors are required to be one (with
the initial set of given vectors of norm one) the authors obtained
a lower bound for the number of unit norm vectors we have to add
(\cite{FWW},Theorem 1.2); but they showed that their lower bound
is not sharp in some cases.

In this note, we calculate the minimum number of vectors we have
to add to $\cF$ to obtain a \emph{tight completion}. Moreover, we
do not require the vectors to be of  norm one; we look for
tight completions with sequences of vectors whose squared norms
are prescribed  by a non-increasing sequence of positive numbers.

Note
that this problem may not have a positive solution for a given set
of initial vectors and a fixed sequence of ``prescribed norms".
Therefore we first find conditions for such a tight frame
completion to  exist. The main tool used here is Theorem \ref{SH}, which relates the squared norms of the vectors in a Bessel sequence with the spectrum of its frame operator.

In order to  state our main results,  we fix some notation used
throughout the paper. Let $\hil$ be a real or complex finite
dimensional vector space with $\dim \hil=n\in \NN$. Let $\f
\subseteq \hil$ be a finite sequence with frame operator $\fop
\cF$ whose eigenvalues (counted with multiplicity) are
$\lambda_1\geq \ldots\geq \lambda_n$, and let $\anorm$ be a
non-increasing sequence of positive real numbers. Finally, let
$\alpha=\tr(\fop \cF)$.
\begin{teo*}[A] Given  $r\in \NN$,  there exists
$\g \subseteq \hil$ such that $\cF\cup\cG$ is a tight frame if and only
if $\frac{1}{n} \left(\suma r + \al \right)\geq
   \la_1$ and
   \begin{equation}
       \frac{1}{n} \left(\suma r + \al \right) \geq \frac{1}{k}\sumal k , \;
        1\leq k \leq \min\{n ,r\}.
\end{equation} On the other hand, there exists an infinite Bessel sequence
$\cG=\{g_i\}_{i=1}^\infty$ in $\hil$ such that $\cF\cup\cG$ is a
tight frame if and only if $\{a_i\}_{i=1}^\infty\in \ell^1(\NN)$,  $\frac{1}{n} \left(\suma \infty +
\al \right)\geq
   \la_1$ and
\begin{equation}
       \frac{1}{n} \left(\suma \infty + \al \right)
       \geq \frac{1}{k}\sumal k , \;
        1\leq k \leq n.
\end{equation}
\end{teo*} So from Theorem A we get necessary and
sufficient conditions for the existence of a sequence $\g$, for
some $r\in \NN\cup\{\infty\}$, with $\|g_i\|^2=a_i$, and such that $\cF\cup\cG$ is a tight frame (for some
suitable constant). If such a completion exists we say that $\cF$
is $(\Aa,r)$-completable. In case $\cF$ is $(\Aa,r)$-completable,
we are then interested in computing the minimum number $r_0$ of
vectors we have to add. In order to state our next result we
introduce the following numbers: let $c_0=\lambda_1$ and for
$1\leq k\leq n$ let
\begin{equation}c_k=\max \left( c_{k-1} ,\;  \frac{1}{k} \sumal k \right).
\end{equation}

\begin{teo*}[B] Assume that $\cF$ is $(\Aa,r)$-completable for
some $r\in \NN\cup\{\infty\}$ and let $r_0\in \NN\cup\{\infty\}$
be the minimum  such that
 $\cF$ is $(\Aa,r_0)$-completable.

 Then
\begin{description}
    \item[Case 1] $r_0<n$ if and only if  $c_{r_0}=\frac{1}{n}\left( \suma {r_0} +\al
    \right)$.
  \item[Case2] $n \leq r_0 <\infty$ if and only if $c_k \neq \frac{1}{n}\left( \suma k +\al \right)\, \forall \, 1\leq k \leq n-1$ and  $r_0$ is the minimum
  such that $c_n \leq \frac{1}{n}\left( \suma {r_0} +\al \right)$.

\item[Case 3] $r_0=\infty$ if and only if  $c_k \neq
\frac{1}{n}\left( \suma k +\al \right)$ for all
  $1\leq k \leq n-1$  and $c_n =\frac{1}{n}\left( \suma
{\infty} +\al \right)$.
\end{description}
\end{teo*}

We should remark that although Theorems A and B are of practical
interest, they are not efficiently (fast) algorithmic
implementable in a computer (see the discussion at the beginning
of Section \ref{algoritmos}). In Section \ref{algoritmos} we deal
with the problem of finding a not so optimal but efficiently
algorithmic computable finite tight completion as follows:
\begin{teo*}[C] Assume that $\Aa$ is a divergent sequence. Let
$d\in \RR$ be an algorithmic computable upper bound for $\|\fop
\cF\|$ and let $c=\max(d+1, d+a_1)$. If $r\in \NN$ is such that $$
\sum_{i=1}^{r-1}a_i< c\cdot n-\tr(\fop \cF) \leq \sum_{i=1}^r a_i
$$ then there exists an algorithmic computable sequence $\g$ such
that $\cF\cup\cG$ is a tight frame and such that $ \|g_i\|^2=a_i$
for $1\leq i\leq r$.
\end{teo*}

We also consider particular cases of Theorems A and B
when $a_i=1 $ for $i\geq 1$.

    \section{Preliminaries on frames and majorization}\label{prelis}

 Throughout the paper,  $\hil$ will be a finite  dimensional (real or complex)
Hilbert space with $\dim \hil=n\in \NN$ and   $\posop$ will denote the cone of bounded positive semi-definite operators on $\hil$. Given  $m\in \NN\cup\{\infty\}$,
a sequence $\cF=\{f_i\}_{i=1}^m \subset \hil$  is a {\sl frame}
for $\hil$ if there exist numbers $ a,b>0$  such that, for every  $f \in \hil$,
\begin{equation}\label{frame}
a\,\|f\|^2\leq \sum_{i=1}^m |\pint{f , f_i}|^2 \leq b\,\|f\|^2
\end{equation}
The optimal  constants in \eqref{frame} are called the \textit{frame bounds}.
If the frame bounds $a,b$ coincide, the frame is called {\sl
$a$-tight} (or simply tight). Finally, tight frames with all its
elements having the same norm are called {\sl equal norm tight
frames}.

The sequence $\cF$ is {\sl Bessel} if there exists $b>0$ such that
the upper bound condition in \eqref{frame} is satisfied. Given a
Bessel sequence $\cF$,  we define its {\sl frame operator}
by
\begin{equation}
\fop \cF f =\sum_{i=1}^m \pint{f , f_i}f_i.
\end{equation}
 It is easy to see that $\fop \cF$ is a positive
semi-definite bounded operator on $\hil$. Moreover, $\cF$ is a
frame if and only if its frame operator $\fop \cF$ is invertible. Indeed, the
optimal frame bounds $a,b$ in \eqref{frame} are respectively
$\la_{\min}(\fop \cF )$ and $\la_{\max }(\fop \cF)$, the minimum
and maximum eigenvalues of $\fop \cF$. In particular, a frame
$\cF$ is $a$-tight if and only if $\fop \cF =aI$. For an
introduction to the theory of frames and related topics see the
books \cite{liChr,{Y}}.

Given a Bessel sequence $\cF$, there is a close relationship
between the norms of its elements and the spectrum of $\fop \cF$
that can be expressed in terms of majorization (see \cite{AMRS} for
details). First, we introduce some definitions. We say that a
sequence $\{a_i\}_{i=1}^m$ is {\sl summable} if  $m\in\NN$,
or if $m=\infty$ and $\{a_i\}_{i=1}^\infty\in \ell^1(\NN)$.

\begin{fed}\label{mayo trucha}
Let $\Aa =\{a_i\}_{i=1}^m$, $\textbf{b}=\{b_i\}_{i=1}^s$ be
 non-increasing summable sequences  of non-negative numbers, with $s,m
\in \NN \cup \{\infty\}$, and let $t=\min\{s,m\}$. We say that
 $\textbf{b}$ majorizes $ \Aa$, noted $ \textbf{b}\succ\Aa $, if
\begin{equation}
\sum_{i=1}^j b_i \geq \suma j  \; \text{ for } 1\leq j\leq t \;
\text{ and } \sum_{i=1}^s b_i=\suma m.
\end{equation}
\end{fed}

If $m=s\in \NN$ in Definition \ref{mayo trucha} then this notion
coincides with the usual vector majorization in $\RR^m$ between
vectors with non-negative entries which are arranged in
non-increasing order (see \cite{horn2}).

%We introduce this notion
%since we want to compare vectors of different sizes in terms of
%(an extended notion of) majorization, avoiding the completion of
%the vectors with entries of zeros (for a detailed analysis of this
%notion of majorization see \cite{{AMRS},{neu}}).

On the other
hand, as an immediate consequence of Definition \ref{mayo trucha}
we see that if $s\in \NN$, and  then $\Aa \prec \b$ if and only if $\Aa
\prec (\b,0_n)$ for every $n\in \NN$, where $(\b,0_n)\in
\RR^{s+n}$, and similarly $(\Aa,0_n)\prec \b$ if $m\in \MM$.

\begin{rem}\label{abajo el chico}
Let $\Aa,\,\textbf{b}$ be as in Definition \ref{mayo trucha},
with $ \textbf{b}\succ \Aa$ and $m<s$. Then $b_i=0$   for
$m+1\leq i\leq s$, since $$\sum_{j=1}^m a_j=\sum_{j=1}^s b_j \geq
b_i + \sum_{j=1}^m b_j \geq \sum_{j=1}^m a_j$$ which implies that
$b_i=0$ since $\b$ has non-negative entries.
\end{rem}

Now we can  state  the frame version of the Schur-Horn theorem,
which we shall need in the sequel.

\begin{teo}\label{SH}
Let $\Aa =\{a_i\}_{i=1}^m$ be a non-increasing sequence of
positive numbers and let $S\in \posop$ with eigenvalues (counted
with multiplicity and arranged in non-increasing order) $\mathbf
\Lambda=\{\la_j\}_{j=1}^n$. Then the following statements are
equivalent:
\begin{enumerate}
    \item  $\Aa \prec \mathbf \Lambda$.
    \item There exists a Bessel sequence $\cG=\{g_i\}_{i=1}^m\subset \hil$
     such that $\|g_i\|^2=a_i$ for $1\leq i\leq m$ and $\fop \cG =
     S$.
\end{enumerate}
\end{teo}

\begin{proof}
If we assume that $S>0$ then the case when $m\in \NN$ is Theorem
4.6 in \cite{AMRS}, while the case when $m=\infty$ is Theorem 4.7 in \cite{AMRS}.
 If the spectrum of $S$ has zeros (note that this is the case
 whenever $m<n$)
 we can reduce
 to the invertible case,
 restricting $S$ to the orthogonal complement of $\ker S$.
\end{proof}

\begin{rem}
We have now a way to look at the problem of tight completions of a
given set of vectors: a set $\cF=\{f_i\}_{i=1}^p $ has a $c$-tight
completion $\cG=\{g_i\}_{i=1}^m$ if and only if $\fop \cG =
cI-\fop \cF$. Thus, by Theorem \ref{SH},this is equivalently to
the fact that the squared norms of $\{g_i\}_{i=1}^m$ are majorized
by the non-increasing sequence $c-\la_n \geq c-\la_{n-1}\geq
\ldots \geq c-\la_1$, where $\la_i$ are the eigenvalues of $\fop
\cF$ (counted with multiplicity and rearranged in decreasing
order).
\end{rem}

\section{Completing a Bessel sequence to a tight frame with prescribed
norms}\label{completing}

%In the rest of the paper $\f$ will denote a finite set of vectors
%in $\hil$  with frame
%operator $\fop \cF$ whose eigenvalues (counted with multiplicity)
%are $\lambda_1\geq \ldots\geq \lambda_n$, $\tr(\fop \cF)=\alpha$,
%and $\Aa=\{a_i\}_{i=1}^\infty$ will denote a non-increasing
%sequence of positive numbers. Recall that $\hil$ is a finite
%dimensional (real or complex) Hilbert space with $\dim\hil=n\in
%\NN$.

\begin{fed}
We say that $\cF$ is {\bf ($\Aa, r$)-completable} if there exists
$r\in \NN\cup\{\infty\}$ and a Bessel sequence $\g\subset \hil$,
with $\|g_i\|^2=a_i$ for $1\leq i\leq r$, and such that $\cF \cup
\cG$ is a tight frame. We say that $\g$ is an {\bf ($\Aa,
r$)-completion} of $\cF$.
\end{fed}

\begin{rem}\label{c}
If $\g$ is an ($\Aa, r$)-completion of $\cF$ then the frame bound
$c\in\RR$ for $\cF\cup\cG$ is determined by the number $r\in
\NN\cup\{\infty\}$ and the norms of the vectors of $\cF$. In fact
$\tr (\fop \gyf )=nc $, and simple computations show that
$$\displaystyle{\tr (\fop \gyf )=\sum_{i=1}^p \|f_i\|^2 + \suma r
}$$ so we have that $c=\frac{1}{n} (\suma r + \tr ( \fop \cF) )$.
In particular, if $r=\infty$ then $\Aa$ is summable.
\end{rem}

 For the sake of clarity in the exposition, in what follows we consider
separately the cases where  \comp r for some $r\in \NN$ and the case
$r=\infty$, although there is no substantial difference in the
arguments involved.

\subsection{Completing with a finite number of vectors}

\begin{teo}\label{caso general}
Let $r\in \NN$. Then \comp r  if and only if
   $\frac{1}{n} \left(\suma r + \al \right)\geq
   \la_1$ and
   \begin{equation}\label{desigualdades}
       %\;\text{ and }\;
       \frac{1}{n} \left(\suma r + \al \right) \geq \frac{1}{k}\sumal k , \;
        1\leq k \leq \min\{n ,r\}.
\end{equation}
\end{teo}

\begin{proof}
Assume that there exists $r\in \NN$ and a finite sequence  $\g$
such that $\fop \gyf=\fop \cF+\fop \cG =cI$ and $\|g_i\|^2=a_i$
for $1\leq i\leq r$. Then $cI-\fop \cF=\fop \cG\geq 0$; in
particular we have $c\geq \|S\|=\lambda_1$. On the other hand, we
see that the eigenvalues of $\fop \cG$ arranged in non-increasing
order are
 $c-\la_n \geq \ldots \geq c-\la_1\geq 0$.
By Theorem \ref{SH} we have \begin{equation}\label{mayo1}(c-\la_n,
c-\la_{n-1}, \ldots , c-\la_1) \succ
(a_1,\ldots,a_r).\end{equation} Then, by Definition \ref{mayo
trucha} we see that $\frac{1}{n} \left(\suma r + \al \right)\geq
   \la_1$ and \eqref{desigualdades} hold, using that
$c=\frac{1}{n} (\suma r + \al )$ by Remark \ref{c}.

Conversely assume that $\frac{1}{n} \left(\suma r + \al
\right)\geq \la_1$ and \eqref{desigualdades}
    hold for $r\in \NN$.  Set  $c
=\frac{1}{n} (\suma {r} + \al )$ and note that the spectrum of the
positive operator $cI-\fop \cF$, $(c-\la_n, c-\la_{n-1}, \ldots ,
c-\la_1)$, majorizes (in the sense of Definition \ref{mayo
trucha}) $\{a_i\}_{i=1}^r$. By Theorem \ref{SH} we conclude that
there exists a finite sequence $\g$ with $\fop \cG=cI-\fop \cF$
and $\|g_i\|^2=a_i$ for $1\leq i\leq r$ and we are done.
\end{proof}

\begin{rem}\label{con finitos}As a consequence of Theorem \ref{caso
general} we see that if $\sum_{i=1}^\infty a_i$ diverges, then
every set of initial vectors $\cF$ is ($\Aa , r$)-completable for
some $r\in \NN$. We shall consider this in section \ref{norma 1}
where we have $a_i=1$ for $i\in \NN$.
\end{rem}

By inspection of the proof of Theorem \ref{caso general} and
Remark \ref{abajo el chico}, we have the following corollaries.
\begin{cor}\label{r<n}
Using the notations of Theorem \ref{caso general},  \comp r
with $r<n$ if and only if, for $ 1\leq i \leq n-r$ and $1\leq k
\leq r$,
\begin{equation}\label{caso 1} \la_i =\frac{1}{n} \left(\suma r  + \al \right),  \, \text{ and }
\,
    \la_1 \geq \frac{1}{k} \sumal k .
\end{equation}
\end{cor}
%\begin{proof}[Proof of Theorem \ref{caso r geq n}]
%As before, the existence of $\g$ satisfying hypothesis 1) is  equivalent (via Theorem %\ref{SH}) to
%$$(c-\la_n,\ldots ,c-\la_1)\succ (a_1,\ldots,a_r)$$

%Note that $c-\la_1\geq 0$ since it is the smallest eigenvalue of $\fop \cG$.
%Again, using definition of majorization, this is equivalent to
%$$kc\geq \suma k +\sumlak \quad \text{ for } \quad 1\leq k \leq n$$
% and
%$$nc-\sumla n= nc-\al =\suma r  $$
% Finally, replacing $c$ by $c=\frac{1}{n} (\suma r + \al )$, we obtain the equivalence of %1) and 2).
%The case $r=\infty$ is identical, with the assumption $\suma \infty <\infty$.

%\end{proof}
%Using Theorem \ref{caso general} and Corollary \ref{r<n}, we get the next corollary, which describes the possibilities for the set $\{ r\in \NN : \text{ \comp r }\}.$

\begin{cor}\label{conjunto de comp}
Let $\cF$ be ($\Aa, r$)- completable for some $r\in \NN$. Then
\begin{enumerate}
\item if $r<n$ then  $\cF$ is not ($\Aa, k$)-completable
for any $k<n$ other than $r$,
\item if $r\geq n$ then \comp k for
every $k\in \NN$ with $ k\geq r$.
\end{enumerate}
\end{cor}

The next result gives different equivalent conditions for a
sequence $\Aa$ and vectors $\cF$ in order to be
$(\mathbf{a},r)$-completable for some $r\in \NN$.
 First, we define
inductively the following numbers: let $c_0=\lambda_1$ and for
$1\leq k\leq n$ let
\begin{equation}\label{los cr} c_k=\max \left( c_{k-1} ,\;  \frac{1}{k}
 \sumal k \right).
\end{equation} It is clear from  definition that
$\lambda_1\leq c_1\leq \ldots\leq c_n$.

\begin{pro}\label{compatibilidad}Let $r\in \NN$. \comp r if and only if
\begin{equation}\label{eqcomp}
\frac{1}{n}\left( \suma r + \al \right)=c_r \text{ for } r<n
\end{equation}
or \begin{equation}\label{eqcomp2} \frac{1}{n}\left( \suma r
+ \al \right)\geq c_n \text{ for } r\geq n.
\end{equation}
Moreover, if $c_r=\frac{1}{n}\left( \suma r + \al \right)$ for
some  $r<n$, then $c_r=\la_1$.
\end{pro}

\begin{proof} Assume that \comp r.
If $r<n$ note that, by \eqref{caso 1} in Corollary \ref{r<n}, we
have $\lambda_1=c_0\leq \ldots\leq c_r=\lambda_1$ and
$\lambda_1=\frac{1}{n}\left( \suma r+ \al \right)$, so
\eqref{eqcomp} holds. If $r\geq n$ then $\min\{n,r\}=n$ and
Theorem \ref{caso general} together with the definition of $c_n$
imply that $$ \frac1n\left(\sum_{i=1}^r a_i+\alpha \right)\geq
c_n.$$ So in this case \eqref{eqcomp2} holds. Conversely, if we
assume \eqref{eqcomp2}, then it is clear \comp r, by Theorem
\ref{caso general}. Assume now that for some $r< n$,
$c_r=\frac{1}{n}\left( \suma r +\al \right )$. We show that \comp
r; indeed, since $nc_r=\suma r +\al$, then
\[ rc_r + (n-r)c_r - \sum_{i=1}^{n-r} \la_i = \suma r + \sumla r \]
so by definition of $c_r$ we have
\[  \sumal r  \leq r \, c_r=\sumal r -  \sum_{i=1}^{n-r}
(c_r-\la_i) \leq \sumal r.\]
But then $$ \la_i=\frac{1}{n}\left(
\suma r +\al \right )\text{ for } 1\leq i\leq n-r $$
 and
 $$\la_1\geq \max_{1\leq k\leq r}\frac{1}{k} \sumal k,$$
 so \comp r, by Corollary \ref{r<n}. The last claim of the proposition is clear from our previous computations.
\end{proof}
We are now able to give a formula for the minimum $r\in \NN$ such
that \comp r, when such an $r\in \NN$ exists.

\begin{teo}\label{optimalidad}
Let $\cF$ be a $(\Aa,r)$-completable  for some $r\in \NN$. Let $r_0\in \NN$ be
the minimum  such that \comp{r_0}. Then
\begin{description}
    \item[Case 1] $r_0<n$ if and only if $c_{r_0}=\frac{1}{n}\left( \suma {r_0} +\al \right)$
  \item[Case2] $r_0\geq n$ if and only if $c_k \neq \frac{1}{n}\left( \suma k +\al \right)$ for all
  $1\leq k \leq n-1$ and $r_0\in \NN$ is the minimum
  such that $c_n \leq \frac{1}{n}\left( \suma {r_0} +\al \right)$.

\end{description}
\end{teo}

\begin{proof} Note that, by Proposition \ref{compatibilidad}, at
least one the cases has to be fulfilled by some $r\in \NN$. If we
assume that case 1 holds for some $r<n$ then, by Proposition
\ref{compatibilidad}, \comp r. By Corollary \ref{conjunto de comp}
case 1 does not hold for $k<n$ with $r\neq k$. It is clear that in
this case $r_0=r$.

Assume now that there is no $r<n$ satisfying case 1 above. Then,
there exists $r\in \NN$ such that $c_n \leq \frac{1}{n}\left(
\suma {r} +\al \right)$; by Proposition \ref{compatibilidad} we
see that \comp r. It is clear that $r_0$ is the minimum natural
number $r$ satisfying this condition. Finally note that if $r\in
\NN$ is such that $c_n \leq \frac{1}{n}( \suma {r} +\al )$ then $$
\frac{1}{n}( \suma {n} +\al ) \leq c_n \leq \frac{1}{n}( \suma {r}
+\al )\ \Rightarrow \suma {n}\leq \suma {r}$$ and $r\geq n$ since
for every $i\in\NN$, $a_i>0$.
\end{proof}

%\begin{exa}[se completa solo con infinitos]
%
%\end{exa}

The next  example  shows that it is possible to obtain a set of
vectors $\cF$ and a sequence $\Aa$ such that \comp r for only one
$r\in \NN$ (in virtue of Corollary \ref{conjunto de comp}, $r<n$).

\begin{exa}
Let $\cF=\{ \sqrt{2}e_1 , \sqrt{2}e_2 , e_3\}$ in $\CC^3$ where
$\{e_i\}$ is the canonical orthonormal basis and let
$\Aa=\{\left(\frac{1}{4}\right)^{i-1}\}_{i=1}^\infty$. Then, easy
computations show that the eigenvalues of $\fop \cF$ are
$\la_1=2$, $\la_2=2$ and $\la_3=1$, so $\al=\tr \fop \cF = 5$. By
Corollary \ref{r<n} \comp 1 since $\la_1=\frac{1}{3}(a_1+\al )$
and $\la_1 \geq a_1 +\la_3$. Moreover, it is clear that if we add
the vector $e_3$ to $\cF$ we obtain a $2$-tight frame.

On the other hand, it easy to see that  $\frac{1}{3}(\suma \infty
+\al )= \frac{19}{9} < \frac{17}{8}= c_3$ so, by Proposition
\ref{compatibilidad}, $\cF$ is not ($\Aa , r $)-completable for
any $r\geq 3$.

\end{exa}

In fact, as the following proposition shows, if $\cF$ is   ($\Aa ,
r $)-completable with $r<n$, the existence of some $r_1\geq
n$ such that $\cF$ is ($\Aa , r_1 $)-completable  depends only on
the tail of the sequence, $\{a_i\}_{i=r+1}^\infty$.

\begin{pro}\label{cola de serie}
Let $\cF$ be $(\Aa,r)$-completable for some $r<n$. There
exists $r_1\in \NN$ with $r_1\geq n$ and such that \comp{r_1} if
and only if $$\frac{1}{n} \sum_{i=r+1}^{r_1} a_i \geq
\max_{r+1\leq k \leq n} \frac{1}{k}\sum_{i=r+1}^k a_i.$$
\end{pro}

\begin{proof}
By Theorem \ref{caso general}, \comp {r_1}  if and only if $$
\frac{1}{n} \left(\suma {r_1} + \al \right)\geq \la_1 \, \text{
and }\, \frac{1}{n} \left(\suma {r_1} + \al \right) \geq
\frac{1}{k}\sumal k  ,\  1\leq k \leq n. $$ By hypothesis and
Corollary \ref{r<n}, $$\la_i =\frac{1}{n} \left(\suma r + \al
\right) ,\ 1\leq i \leq n-r  \, \text{ and } \,
    \la_1 \geq \frac{1}{k}\sumal k  ,\  1\leq k \leq r $$
since \comp r with $r<n$. So \comp {r_1} if and only if $$\quad
\frac{1}{n} \left(\suma {r_1} + \al \right) \geq \frac{1}{k}\sumal
k  ,\ r+1\leq k \leq n$$ or equivalently, if for every
$r+1\leq k \leq n$
\begin{align*}
%\frac{1}{n} \left(\suma {r_1} + \al \right) & \geq \frac{1}{k}\sumal k  \\
 \suma r + \al + \sum_{i=r+1}^{r_1}a_i &
\geq \frac{n}{k}\left(\suma r + \sum_{i=r+1}^{k}a_i +\alpha
-(n-k)\la_1\right)\\  \suma r + \al + \sum_{i=r+1}^{r_1}a_i & \geq
\frac{n}{k}\left(\suma r +\al \right)+
\frac{n}{k}\sum_{i=r+1}^{k}a_i -\frac{n-k}{k}\left(\suma r + \al
\right)\\  \suma r + \al + \sum_{i=r+1}^{r_1}a_i & \geq \suma r +
\al +\frac{n}{k}\sum_{i=r+1}^{k}a_i \\ \sum_{i=r+1}^{r_1}a_i &
\geq \frac{n}{k}\sum_{i=r+1}^{k}a_i \ ,
\end{align*}
since by hypothesis $\la_i =\frac{1}{n} (\suma r + \al)$  for
$1\leq i \leq n-r$.
\end{proof}

\subsection{Completing with a infinite number of vectors.
Proof of Theorems A and B}

In this section we consider some complementary results to those
obtained in the previous section and prove  Theorems A and B.

If $\f$ and $\Aa$ are as before, then a necessary condition for
$\cF$ to be $(\Aa,\infty)$-completable is that $\Aa\in
\ell^1(\NN)$ (Remark \ref{c}).
\begin{teo}\label{caso general infty} \comp \infty  $ \ $(by a Bessel sequence) if and
only if $\Aa\in \ell^1(\NN)$, $\frac{1}{n} \left(\suma \infty +
\al \right)\geq
   \la_1$ and
   \begin{equation}\label{desigualdades infty}
       %\;\text{ and }\;
       \frac{1}{n} \left(\suma \infty + \al \right) \geq \frac{1}{k}\sumal k , \;
        1\leq k \leq n,
\end{equation} or equivalently if $(\Aa\in \ell^1(\NN))$
\begin{equation}\label{eqcomp2 infty} \frac{1}{n}\left( \suma
\infty + \al \right)\geq c_n.
\end{equation}
\end{teo}
The proof of Theorem \ref{caso general infty}, which is based on
Theorem \ref{SH}, is similar to that of Theorem \ref{caso general}
and Proposition \ref{compatibilidad}; we leave the details to the
interested reader.

\medskip

\noindent {\it Proof of Theorem (A)}. The first part of the
theorem is Theorem \ref{caso general}, while the second part is
Theorem \ref{caso general infty}. \qed

\medskip

\noindent {\it Proof of Theorem (B)}. Assume there exists a
natural number $r\in \NN$ such that \comp r. Then $r_0\leq r$ and
in this case the theorem follows from Theorem \ref{optimalidad}.
If there is no $r\in \NN$ such that \comp r, then \comp \infty $\
$so by Theorem \ref{caso general infty} $\Aa\in \ell^1(\NN)$ and
$\frac{1}{n}( \suma \infty + \al )\geq c_n$. If $\frac{1}{n}(
\suma \infty + \al )>c_n$ then there exists $r\in \NN$ such that
$\frac{1}{n}( \suma r + \al )\geq c_n$ but then, by Proposition
\ref{optimalidad} we get that \comp r,  a contradiction.
\qed

We finish with the counter-part of Proposition \ref{cola de serie}
for the infinite completion case.

\begin{pro}\label{cola de serie infty}
Let $\Aa\in \ell^1(\NN)$ and let $\cF$ be $(\Aa,r)$-completable for
some $r<n$. Then, \comp{ \infty} if and only if $$\frac{1}{n}
\sum_{i=r+1}^{\infty} a_i \geq \max_{r+1\leq k \leq n}
\frac{1}{k}\sum_{i=r+1}^k a_i.$$
\end{pro}

\section{Equal norm tight frames}\label{norma 1}

In this section we consider the particular case when $\anorm$ is a
constant sequence, $a_i=1$ for all $i\in \NN$ (the general case
follows in an analogous way). Note that in this case
 \comp r for some $r\in \NN$; so we shall compute the minimum natural
number $r$ of vectors with norm one we have to add to $\cF$ in
order to get a tight frame. We keep the notation of the previous
section for $\f$, $\fop \cF$, $\lambda_1\geq \ldots\geq \lambda_n$
and $\alpha$.

%As before, in this section $\f$ denotes a finite sequence of
%vectors in $\hil$ not necessarily of norm one, with frame operator
%$\fop \cF$ whose eigenvalues (counted with multiplicity) are
%$\lambda_1\geq \ldots\geq \lambda_n$. Moreover, let $\tr(\fop
%\cF)=\al$.

\begin{rem}\label{los cr en uno}
Under our present assumption that $a_i=1$ for all $i\in \NN$ we
have that $$c_r=\max\left( \lambda_1, 1+\frac{1}{r}\sumla {r}
\right). $$ Indeed, if $j\leq k$ then $\frac{1}{j}\sumla j \leq
\frac{1}{k}\sumla k$  since $\la_1\geq \ldots \geq \la_n\geq 0$.
\end{rem}

This simpler formula for the coefficients $c_j$ provides the
following characterization for the optimal number of elements for
tight completions with norm one vectors.

%The following result is the re-statement of Theorem
%\ref{optimalidad} using Remark \ref{los cr en uno}.

\begin{teo}\label{el minimo}
Let $h:=\sum_{i=2}^n \la_1-\la_i$, and denote by $r_0$
 the minimum number of norm one vectors we have to add to
$\cF$ in order to have a tight frame.
\begin{description}
    \item[Case 1] Suppose that $h<n$. Then $r_0=h\, $ if   $ \, h\in \NN$  and
    $1+\frac{1}{h}\sum_{i=1}^h\la_{n-i+1}\leq \la_1$ (in particular,
    $c_h=\la_1$).  Otherwise, $r_0=n$.
\item[Case 2] If $h\geq n$ , $r_0$ is the minimum integer greater
than or equal to $h$. \end{description}
\end{teo}

\begin{proof}
Assume that  $h<n$; then, since $h=n\la_1-\al$, we have that
$c_n=1+\frac{\al}{n}$ by Remark \ref{los cr en uno}. If in
addition  $h\in \NN$ and $1+\frac{1}{h}\sum_{i=1}^h \la_{n-i+1}
\leq \la_1$, so $c_h=\frac{1}{n}(h+\al)=\la_1$, then $r_0=h$ by
Theorem \ref{optimalidad}. Otherwise, $c_k\neq \frac{1}{n}(k+\al)$
for all $k<n$ (if $c_k=\frac{1}{n}(k+\al)$ for some $k<n$, then by
Proposition \ref{compatibilidad} $c_k=\la_1$ and $h$ would be a
natural number); since $c_n=1+\frac{\al}{n}$, the minimum integer
greater than or equal to $nc_n-\al$ is $n$ so $r_0=n$  by  Theorem
\ref{optimalidad}.

Finally,  $h\geq n$   implies that $c_k\neq \frac{1}{n}(k+\al)$
for all $k<n$ and $c_n=\la_1$. Therefore, again by Theorem
\ref{optimalidad}, $r_0$ is the minimum integer greater than or
equal to $n\la_1-\al=h$.
\end{proof}

\begin{rem}
%Note that the number $h=\sum_{i=2}^n \la_1-\la_i$ in some sense measures how ``far of being tight" is the set of %vectors $\cF$. Theorem
Note that, as a consequence of Theorem \ref{el minimo}, if $\la_1-\la_2\geq \frac{n}{n-1}$ then $h=\sum_{i=2}^n \la_1-\la_i\geq n$, so $\cF$ can not be completed to a tight frame with less than $n$ unit norm vectors.

In addition, the number $h$ can be seen as a kind of measure of how ``far of being tight" is the set of vectors $\cF$, in the sense that $h=0$ if and only if the set $\cF$ is already a tight frame.
\end{rem}

\begin{exa}
This example is taken from  \cite{FWW}. It is interesting because
it shows the difference between the cases when we can complete
$\cF$ to a tight frame with $r<n$ or $r\geq n$ vectors. Let
$f_1=(1,0)$ and $f_2=(\cos \theta , \sin \theta)$ in $\RR^2$, and
consider $a_i=1 \; \forall i$. It easy to see that the eigenvalues
of $\fop \cF$ are $1\pm \cos \theta$, hence
$h=\la_1-\la_2=2\,|\cos \theta \,|$. Therefore, by Theorem \ref{el
minimo}, the minimum number $r_0$ of unit vectors we have to add
to obtain a tight frame is 2, unless  $\theta=\frac{2}{3}\pi$ or
$\theta=\frac{4}{3}\pi$, where $r_0=1$. Note that when $r_0=1$ the
tight frame obtained is the well known ``Mercedes Benz" (it is
--up to rigid rotations, reflections and negation of individual
vectors-- the only unit norm tight frame on $\RR^2$ with three
elements \cite{GKK}).

\end{exa}

A consequence of Theorem \ref{el minimo} is the characterization of  the minimum
number of vectors that we have to add in order to get a tight
frame, in the particular case when $\cF$ is a unit norm tight
frame on its linear span.

\begin{pro} Let $\f$ be a
 unit norm $\frac{p}{d}$-tight
 frame on its span, where $d<n$ is the
 dimension of span $\cF$. Then, the minimum number
 $r_0$ of unit norm vectors we have
 to add to $\cF$ in order to obtain a tight frame in $\hil$ is:
\begin{itemize}
    \item[a)] $(n-d)\frac{p}{d}\; \text { if } \; (n-d)\frac{p}{d}< n$ and
    $(n-d)\frac{p}{d}\in \NN $.
  \item[b)] $n \; \text{ if } \; (n-d)\frac{p}{d}< n$
  and $(n-d)\frac{p}{d} \notin \mathbb{N} $.
  \item[c)] the minimum integer
  greater than or equal to $(n-d)\frac{p}{d} \; \text{ if } \; (n-d)\frac{p}{d}\geq n$.
\end{itemize}
\end{pro}

\begin{proof}
Since $\cF$ is an  unit norm tight frame on a subspace of
dimension $d$, the eigenvalues of $\fop \cF$ are:
$\la_i=\frac{p}{d}\geq 1$ for   $ 1\leq i\leq d$, and $\la_i=0$
for $d+1\leq i\leq n$.  Therefore, $h=\sum_{i=2}^n\la_1-\la_i = (n-d)\frac{p}{d}$. Moreover, if $h<n$ and $h\in \NN$,  then  $1+\frac{1}{h}\sum_{i=1}^h\la_{n-i+1}= \la_1$. Indeed,
\begin{equation}\label{ecuaX25}1+\frac{1}{h}\sumla h=
1+\frac{h-(n-d)}{h}\frac{p}{d}=\frac{p}{d}\end{equation}
the proposition is then a consequence of Theorem \ref{el minimo}.
\end{proof}

\section{Some remarks regarding algorithmic
implementation}\label{algoritmos}

Let $\f\subseteq \hil$ and assume that $\Aa$ is a divergent
sequence. Then, by Remark \ref{con finitos}, \comp r for some
$r\in \NN$. From the proof of Theorem \ref{caso general} we see
that, if $c=\frac{1}{n} \left(\suma r + \al \right)$ then equation
\eqref{mayo1} holds. Therefore, by Theorem \ref{SH},
\emph{theoretically}, there exists a Bessel sequence
$\cG=\{g_i\}_{i=1}^r\subset \hil$
     such that $\|g_i\|^2=a_i$ for $1\leq i\leq r$ and $\fop \cG =
     cI-\fop \cF$. In this case, $\cG$ is a $(\Aa,r)$-completion of $\cF$; moreover, if $r\in \NN$ is obtained as in Theorem \ref{optimalidad} then $\cG$ would be $(\Aa,r)$-tight completion having the minimum number of vectors
     for which a tight completion of $\cF$ exists.
     Although constructive, the proof of Theorem \ref{SH}
     is not practicable; it depends on some
     matrix decompositions which can not be performed efficiently by
     a computer for large values of $t=\min\{n,r\}$.

     There are several recent papers related to algorithmic construction of frames with additional properties. In \cite{CL} Casazza and Leon considered the problem of constructing frames with prescribed properties from an algorithmic point of view; in particular, they obtained an algorithm for constructing tight frames with prescribed norms of its elements, under the admissibility conditions of Theorem \ref{SH}.
In \cite{FWW} there is a fast algorithm for constructing tight
frames with prescribed norms of its elements based on Householder
transformations; in \cite{TDHS2}   a fast
algorithmic proof of some results related to the Schur-Horn
theorem is considered and as a consequence a generalized one-sided Bendel-Mickey
algorithm (see Theorem \ref{algo1} below) is obtained. Still, as
far as we know, the problem of constructing a frame for $\hil$
with prescribed general (positive definite) frame operator and
norms (that are admissible in the sense of Theorem \ref{SH}) using
an \emph{efficient} computable algorithm has not been solved: we
remark that for the purposes of this discussion, the
diagonalization of a positive semi-definite matrix is considered
as \emph{not} efficiently computable. If such an algorithm is
obtained, then optimal tight frame completions can be constructed
as described in the first paragraph of this section. In what
follows we shall consider a not so optimal tight frame completion
of a given set $\f$ but that is efficiently algorithmic
computable, based on the generalized one-sided Bendel-Mickey
algorithm and the Cholesky's decomposition.

Let us begin with the
 following result from \cite{TDHS2}. We remark
 that our notation is opposite to that in \cite{TDHS2}
 so we translate their result into our terminology.

\begin{teo}[\cite{TDHS2}]\label{algo1}
Let $\Aa=\{a_i\}_{i=1}^r,\,\Bb=\{b_i\}_{i=1}^r$ be two finite and
non-increasing sequences of positive numbers such that $\Aa\prec
\Bb$. Let $X$ be an $n\times r$ matrix whose squared
columns norms are listed by $\Bb$. Then there is a finite sequence
of algorithmic computable plane rotations $U_1,\ldots,U_{r-1}\in
\M_r(\CC)$ such that $X(U_1\cdots U_{r-1})$ has squared
columns norms listed by $\Aa$.
\end{teo} Actually,
each plane rotation that appears in the theorem above operates
non-trivially in the coordinate plane span$\{e_i,e_j\}$ for some
$1\leq i,\,j\leq r$ (see \cite{TDHS2} for details). Note that the
initial matrix $X$ and the final matrix $Y=
X(U_1\cdots U_{r-1})$ satisfy $X X^*=
YY^*$.

Taking into account Theorem \ref{algo1}, an strategy to construct a
frame with prescribed frame operator $S\in \M_n(\CC)$ and norms of
its elements listed by $\Aa$ (satisfying the conditions in Theorem
\ref{SH}) would be the following: consider a diagonalization
$S=U\text{diag}(\lambda_1,\ldots,\lambda_n)U^*$ and the
factorization $XX^*=S$ with
$X=U\text{diag}(\sqrt{\lambda_1},\ldots, \sqrt{\lambda_n})$. Note
that the squared norms of the columns of $X$ are listed by
$(\lambda_1,\ldots,\lambda_n)$ so we can apply Theorem \ref{algo1}
and obtain $Y=X(U_1\cdots U_{r-1})$ with $ YY^*=S$ with the
squared norms of the columns of $ Y$ given by $\Aa$.
Unfortunately, we consider this procedure as not an efficiently
computable one, so we have to find an alternative approach.

\begin{rem}
In what follows we shall make use of the well known Cholesky's decomposition $S=RR^*$ for a positive definite matrix $S$. Note that in this case Cholesky's decomposition is unique, and there are several strategies for calculating the matrix $R$ is an efficient way.
\end{rem}

\subsection{Algorithm for constructing tight completions}
Along this section we prove Theorem C; we begin with an informal
discussion of the algorithm.

Assume that the non-increasing sequence of positive numbers
$\{a_i\}_{i=1}^\infty$ forms a divergent series, so that \comp t
for some $t\in \NN$. Let $S=\fop \cF$ and let $c>\| S\|$
that we shall consider as a variable. We will obtain an
algorithmic computable value of $c$ for which the Cholesky's
decomposition $cI_n-S=RR^*$ satisfies that the
squared norms of the columns of $R$ mayorize $\{a_i\}_{i=1}^{r}$
for an integer  $r\geq n$. Once we have obtained such  $c$,  we
apply Theorem \ref{algo1} and get a finite sequence $\g$ with
frame operator $cI_n-S$ and $\|g_i\|^2=a_i$, for $1\leq i\leq
r$.

Let $c\geq \|S\|+\beta$  so $\lambda_{\min}(cI-S)=c-\|S\|\geq \beta$, where $\beta>0$ is a fixed number controlling the invertibility of $cI-S$. Let
$R=R(c)$ be the upper triangular matrix obtained from the
Cholesky's decomposition of $cI-S$ (note that the hypothesis on $c$ is made in order that the Cholesky's algorithm becames stable). Then $RR^*=cI-S$ and note
that $c-\|S\|=\lambda_{\min}(RR^*)=\lambda_{\min}(R^*R)$ so, if $C_i(R)$ denotes the $i$-th column of $R$ then
$$\min_{1\leq i\leq n} \|C_i(R)\|^2\geq c-\|S\|,$$
since $\|C_i(R)\|^2=(R^*R)_{ii}$ and $(R^*R)_{ii}\geq \lambda_{\min}(R^*R)$ for $1\leq i\leq n$.
%the eigenvalues (counted with multiplicities) of  $R^*R\in
%M_n(\CC)$ arranged in non-increasing order are $(c-\lambda_n,
%\ldots, c-\|S\|)$, where $(\lambda_1,\ldots,\lambda_n)$ are the
%%eigenvalues (counted with multiplicities) of $S$ arranged in
%non-increasing order.
%S, and
%$((R^*R)_{ii})_{i=1}^n\prec (c-\lambda_i)_{i=1}^n$ then we get
In particular
$\sum_{i=1}^k \|C_i(R)\|^2=\sum_{i=1}^k (R^*R)_{ii}\geq k\cdot
(c-\|S\|)$. Let $c\geq \max(\|S\|+\beta,\|S\|+a_1)$ and note that then
\begin{equation}\label{hola} c\geq \frac{1}{k}\sum_{i=1}^k a_i+\|S\|, \
\ \text{for} \ 1\leq k\leq n\end{equation} since
$\frac{1}{k}\sum_{i=1}^k a_i\geq\frac{1}{h}\sum_{i=1}^h a_i$ if
$1\leq k\leq h\leq n$. Let $r\in \NN$ be such that
\begin{equation}\label{agua}\sum_{i=1}^{r-1} a_i<
\sum_{i=1}^n\|C_i(R(c))\|^2=c\cdot n-\tr(S)\leq
\sum_{i=1}^ra_i\end{equation} so $r\geq n$.
We define $c'=\frac{1}{n}(\sum_{i=1}^r a_i+\tr(\fop \cF))$, where $r$ is defined by \eqref{agua} so that, if $R(c')$ denotes the Cholesky's decomposition of $c'I-\fop \cF$ then
we get
$(a_i)_{i=1}^r\prec (\|C_i(R(c'))\|^2)_{i=1}^n$.

Thus, with this
$c'\in \RR$ and $r\in \NN$ we can apply Theorem \ref{algo1} to the
matrix $X=[R(c'),\ 0_{n\times (r-n)}]$ and get the (efficiently
algorithmic computable) $n\times r$ matrix $Y$ such that $YY^*=S$
and $\|C_i(Y)\|^2=a_i$ for $1\leq i\leq r$; setting $g_i=C_i(Y)$
we get $\{g_i\}_{i=1}^r $ with the desired properties. We briefly
resume the previous considerations in the following pseudo-code
implementation:

\begin{itemize}
\item Find an algorithmic computable upper bound $d$ for $\|S\|$.
\item Compute $c=\max(d+\beta,d+a_1)$ (where $\beta>0$ is previously fixed)
and $r\in \NN$ satisfying \eqref{agua}.
\item Redefine $c:=\frac{1}{n}(\sum_{i=1}^r a_i+\tr(\fop \cF))$.
\item Compute the Cholesky's
decomposition $cI-S=RR^*$.
\item Apply Theorem \ref{algo1} to the $n\times r$ matrix $[R, \
0_{n\times (r-n)}]$ and get the $n\times r$ matrix $Y$ such that
$cI-S=YY^*$ and $\|C_i(Y)\|^2=a_i$ for $1\leq i\leq r$.
\item Define $g_i=C_i(Y)$ for $1\leq i\leq r$.
\end{itemize}

\begin{exa}
Assume that $\|f_i\|=1$ for $1\leq i\leq p$ and that $\|a_i\|=1$, so we are looking for unit norm tight completions of a unit norm family of vectors $\cF$. In this case,  it is shown in \cite{FWW} that if $d=\ent{\|\fop \cF\|+1}$ , where $\ent{h}$ denotes the smallest integer greater than or equal to $h$,  there always exists a unit norm tight completion of $\cF$ with $dn-p$ elements. Our arguments above show that there exists an efficiently algorithmic computable unit norm tight completion with $\ent{ n\cdot(\|\fop \cF\|+1)-p}$ (assuming that we can compute efficiently $\|\fop \cF\|$ and seting $\beta=1$). Note that in general we have that $n\cdot\ent{\|\fop \cF\|+1}-p\geq \ent{ n\cdot(\|\fop \cF\|+1)-p}$ .
\end{exa}

\medskip

\noindent{\bf Acknowledgements.} We would like to thank Professors
Demetrio Stojanoff and Nélida Etchebest for interesting
suggestions regarding the material in this note.

\end{document}